\documentclass{article}
\usepackage{amssymb}

\begin{document}

\title {\large\bf Easily Testable Necessary and Sufficient Algebraic Criteria
for Delay-independent Stability of a Class of Neutral
 Differential Systems
\thanks{Supported by National Natural Science Foundation of China(No. 60204006, No.60334020 and No.
10372002), the Knowledge Innovation Pilot Program of Chinese
Academy of Sciences, the Program of National Key Laboratory of
Intelligent Technology and Systems of Tsinghua University, and
National Key Basic Research and Development Program(No.
2002CB312200). Corresponding author: Long Wang, E-mail:
longwang@mech.pku.edu.cn}}

\author{\normalsize Ping Wei$^{1)}$\quad\quad Qiang Guan$^{1)}$\quad\quad Wensheng Yu$^{1),3)}$\quad\quad Long Wang$^{2)}$ \\[0.2cm]
\small  1)Laboratory of Complex Systems and Intelligence Science, Institute of Automation,\\
\small  Chinese Academy of Sciences, Beijing 100080, P. R. China.\\[0.05cm]
\small  2)Center for Systems and Control, Department of Mechanics and Engineering Science,\\
\small  Peking University, Beijing 100871, P. R. China.\\[0.05cm]
\small  3)National Key Laboratory of Intelligent Technology and Systems, Tsinghua University,\\
\small  Beijing 100084, P. R. China.\\[0.2cm]
}
\date{} \maketitle

\begin{quote} \begin{small} \noindent \footnotesize
{\bf Abstract:} This paper analyzes the eigenvalue distribution of
neutral differential systems and the corresponding difference
systems, and establishes the relationship between the eigenvalue
distribution and delay-independent stability of neutral
differential systems. By using the ``Complete Discrimination
System for Polynomials", easily testable necessary and sufficient
algebraic criteria for delay-independent stability of a class of
neutral differential systems are established. The algebraic
criteria generalize and unify the relevant results in the
literature. Moreover, the maximal delay bound guaranteeing
stability can be determined if the systems are not
delay-independent stable. Some numerical examples are provided to
illustrate the effectiveness of our results.
\end{small}
\end{quote}

\vskip 0.8cm \section{ Introduction} \vskip 0.5cm

Neutral differential system model can often be found in control
process, physics, chemical engineering and ecology
\cite{BC63,Hal77,HV93,QLWZ89}. Stability of neutral differential
systems is an important performance index and has been
investigated by many authors. From the control-theoretic
viewpoint, the delay-independent stability of a delay system
ensures robustness and reliability of the system for all delays.
It is well known that the asymptotic behavior of the zero solution
of a linear autonomous retarded differential system is determined
by its eigenvalues. Just as ordinary differential systems, the
stability of the zero solution is equivalent to all eigenvalues
being with negative real parts
\cite{BC63,GTY01,Hal77,HIT85,HV93,KTO03,QLWZ89,YW00}. However, it
is more complex when considering the neutral differential system,
i.e., all eigenvalues being with negative real parts are not
sufficient for the stability of the zero solution
\cite{BC63,Bru70,Hal77,HIT85,HV93,KTO03,QLWZ89}, and the global
hyperbolicity of a corresponding difference system must be taken
into account \cite{AH80,HIT85,HV93}. In neutral differential
system, when characterizing system stability by means of the
corresponding characteristic equation, it is usually required that
the supremum of the real parts of all eigenvalues is negative
\cite{BC63,Hal77,HIT85,HV93,KTO03,QLWZ89}, which is notably
different from retarded differential system.

Consider the linear autonomous neutral differential system

\begin{equation}
\dot{x}(t)-\sum^{N}_{k=1}B_{k}\dot{x}(t-\gamma_{k}\cdot r)
 =A_{0}x(t)+\sum^{N}_{k=1}A_{k}x(t-\gamma_{k}\cdot r) \label{eq1}
\end{equation}
where $R^+=[0,\infty), R=(-\infty,+\infty), x\in R^{n}, A_{0}\in
R^{n\times n}, A_{k}\in R^{n\times n}, B_{k}\in R^{n\times n},
r=(r_{1},\cdots,r_{M})\in (R^{+})^{M},
\gamma_{k}=(\gamma_{k_{1}},\cdots,\gamma_{k_{M}}),
\gamma_{k_{j}}\geq 0$ are integers, $\gamma_{k}\neq 0,
\gamma_{k}\cdot r=\sum_{j=1}^{M}\gamma_{k_{j}}r_{j},
k=1,2,\cdots,N, j=1,2,\cdots,M$. Suppose
$\det[I-\sum^N_{k=1}B_k]\neq 0$ ($\hbox{the notation} \det[\cdot]$
stands for the determinant of the corresponding matrix). The
characteristic equation of System (1) is
\begin{equation}
G(\lambda,r,A,B)=\det [\lambda (I-\sum^{N}_{k=1}B_{k}e^{-\lambda
\gamma_{k}\cdot r}) -A_{0}-\sum^{N}_{k=1}A_{k}e^{-\lambda
\gamma_{k}\cdot r}]=0 \label{eq2}
\end{equation}
where $I$ stands for the $n\times n$ identity matrix,
$A=(A_0,A_1,\cdots,A_N),B=(B_1,\cdots,B_N)$.

{\bf Definition 1} \cite{HIT85} \ \ System $(\ref{eq1})$ is said
to be {\bf delay-independent stable}, if $\forall r\in (R^{+})^M,$
the supremum of the real part of $\lambda$ satisfying Equation
$(\ref{eq2})$ is negative, that is, there exists a $\delta
>0$ , such that $\{\mbox{Re}\lambda:G(\lambda,r,A,B)=0\}\cap
[-\delta,\infty)=\phi$.

In the literature, delay-independent stability is also called
all-delay stability, unconditional stability or absolute
stability.

Consider the following difference system which is closely related
to System $(\ref{eq1})$
\begin{equation}
x(t)-\sum^{N}_{k=1}B_{k}x(t-\gamma_{k}\cdot r)=0 \label{eq3}
\end{equation}
The characteristic equation of System $(\ref{eq3})$ is

\begin{equation}
E(\lambda,r,B)=\det [I-\sum^{N}_{k=1}B_{k}e^{-\lambda \gamma_{k}\cdot r}]=0 \label{eq4}
\end{equation}

{\bf Definition 2} \cite{HIT85} \ \ System $(\ref{eq3})$ is said
to be {\bf globally hyperbolic} at $B$, if $\forall r\in
(R^{+})^M,$ the real part of $\lambda$ satisfying Equation
$(\ref{eq4})$ does not intersect with a neighborhood of zero, that
is, there exists a $\delta
>0,$ such that $\{\mbox{Re}\lambda:E(\lambda,r,B)=0\}\cap
[-\delta,\delta]=\phi$.

In the papers \cite{HIT85,HV93,KTO03}, some elegant theoretical
results on the stability of the neutral differential system
$(\ref{eq1})$ have been established. One important result is

{\bf Lemma 1} \cite{HIT85} \ \ Let
\begin{equation}
p(\lambda,s_{1},\cdots,s_{M},A,B)
=\det [\lambda I-(I-\sum^{N}_{k=1}B_{k}s^{\gamma_{k_{1}}}_{1} \cdots s^{\gamma_{k_{M}}}_{M})^{-1}
(A_{0}+\sum^{N}_{k=1}A_{k}s^{\gamma_{k_{1}}}_{1} \cdots s^{\gamma_{k_{M}}}_{M})] \label{eq5}
\end{equation}
then, System $(\ref{eq1})$ is delay-independent stability if and
only if

$(H1)$ System $y(t)-\sum^{N}_{k=1}B_{k}y(t-\gamma_{k}\cdot r)=0$
is globally hyperbolic at $B$;

$(H2)$
$\mbox{Re}\lambda[(I-\sum^{N}_{k=1}B_{k})^{-1}\sum^{N}_{k=0}A_{k}]<0$,
where $\mbox{Re}\lambda[\cdot]$ denotes the real part of
eigenvalues of the corresponding matrix;

$(H3)$ $p(iy,s_{1},\cdots,s_{M},A,B)\neq 0,\forall y \in R- \{ 0\}, |s_{j}|=1, j=1,2,\cdots,M$ .

The above lemma gives necessary and sufficient conditions on
delay-independent stability of a general linear autonomous neutral
differential system, which is the fundamental basis for studying
neutral differential systems. However, it is difficult to apply
these conditions in practice, because $(H1)$ and $(H3)$ are both
transcendental and can't be explicitly tested in engineering
computations. Therefore, many authors in the area of differential
systems tried to derive simple easily-testable stability criteria
from various viewpoints \cite{Han01,HH96,
KY88,Li88,PW00,RL02,RZ99,VN97,WLB02}. These papers considered some
special systems and just obtained some sufficient conditions,
furthermore, most of which have never considered the global
hyperbolicity of the corresponding difference systems.

In this paper, we consider a more general linear multi-delay
neutral differential system as follows

\begin{equation}
\dot{x}(t)-\sum^{N}_{k=1}B_{k}\dot{x}(t-k\tau)=A_{0}x(t)+\sum^{N}_{k=1}A_{k}x(t-k\tau) \label{eq6}
\end{equation}
where $x\in R^{n}, A_{0}\in R^{n\times n}, A_{k}\in R^{n\times n},
B_{k}\in R^{n\times n}, \tau \in R^{+}, k=1,2,\cdots,N$, and we
assume $\det[I-\sum^N_{k=1}B_k]\neq 0 $. Based on Lemma 1, and
discussions on the global hyperbolicity of the corresponding
difference system \cite{AH80,HIT85}, algebraic criteria for
delay-independent stability of the above system will be
established by means of the ``Complete Discrimination System for
Polynomials" \cite{Yang99,YHZ96,YZH96}. The algebraic criteria
generalize and unify the relevant results in the literature.
Moreover, the maximal delay bound guaranteeing stability can also
be determined when the system is not delay-independent stable.

This paper is arranged as follows: In section 2, the relevant
results of the ``Complete Discrimination System for Polynomials"
are presented; Section 3 discusses the global hyperbolicity of
difference system $y(t)-\sum^{N}_{k=1}B_{k}y(t-k \tau)=0$; Section
4 analyzes the delay-independent stability of System $(\ref{eq6})$
and establishes the easily testable algebraic criteria; Section 5
presents some corollaries and examples by means of the algebraic
criteria provided in section 4; Finally, a simple conclusion is
given in section 6.

\vskip 0.7cm
\section{Complete Discrimination System for
Polynomials}

From Lemma 1 in section 1, the problem of delay-independent
stability of System $(\ref{eq6})$, boils down to analyzing the
eigenvalues' distribution of System $(\ref{eq6})$ and the
corresponding difference system. In the sequel, we will transform
this problem into the real root distribution of some related
polynomials. The following is some results of the ``Complete
Discrimination System for Polynomials", which can effectively
discriminate the distribution of the roots of polynomials by a
series of explicit expressions and calculations with polynomial
coefficients \cite{Gan59,Yang99,YHZ96,YZH96}.

Suppose
\begin{equation}
f(x)=a_0x^n+a_1x^{n-1}+\cdots +a_n
\end{equation}
 the Sylvester matrix \cite{Gan59,Yang99,YHZ96,YZH96} of $f(x)$ and its derivative
$f^{\prime }(x)$ is defined as
 \begin{center}
$\left[
\begin{array}{cccccccc}
a_0 & a_1 & a_2 & \cdots  & a_n &  &  &  \\
0 & na_0 & (n-1)a_1 & \cdots  & a_{n-1} &  &  &  \\
& a_0 & a_1 & \cdots  & a_{n-1} & a_n &  &  \\
& 0 & na_0 & \cdots  & 2a_{n-2} & a_{n-1} &  &  \\
&  &  & \cdots  & \cdots  &  &  &  \\
&  &  & \cdots  & \cdots  &  &  &  \\
&  &  &  & a_0 & a_1 & \cdots  & a_n \\
&  &  &  & 0 & na_0 & \cdots  & a_{n-1}
\end{array}
\right] $
\end{center}
which is called the {\bf discrimination matrix} of $f(x)$, denoted as $Discr(f)$.

\begin{center}
$[D_1(f),D_2(f),\cdots ,D_n(f)]$
\end{center}
is a sequence of the determinants of the first $n$ even-order
principal sub-matrixes of $Discr(f)$, formed by the first $2k$
rows and first $2k$ columns, $k=1, 2, ..., n$, which is called the
{\bf discrimination sequence} of $f(x)$. Furthermore

\begin{center}
$[sign(D_1),sign(D_2),\cdots ,sign(D_n)]$
\end{center}
is called the {\bf sign list} of the discrimination
sequence$[D_1,D_2,\cdots ,D_n]$, where $sign(\cdot)$ is the sign
function, namely

$$
sign(x)=\left\{
\begin{array}{ll}
1\quad\quad\quad &\mbox{if } x>0,\\
0&\mbox{if } x=0,\\
-1&\mbox{if } x<0 \hbox{.}\\
\end{array}
\right .
$$

Construct a {\bf revised sign list} $[\varepsilon _1,\varepsilon
_2,\cdots ,\varepsilon _n]$ for a given sign list $[s_1,s_2,\cdots
,s_n]$ as follows:

$1)$ If $[s_i,s_{i+1},\cdots ,s_{i+j}]$ is a section of the given
list, where $s_i\neq 0;s_{i+1}=s_{i+2}=\cdots
=s_{i+j-1}=0;s_{i+j}\neq 0,$ then replace the subsection
$$[s_{i+1},s_{i+2},\cdots ,s_{i+j-1}]$$
by
$$[-s_i,-s_i,s_i,s_i,-s_i,-s_i,s_i,s_i,-s_i,\cdots ]\hbox{,}$$
namely, $\varepsilon _{i+r}=(-1)^{[\displaystyle\frac
{r+1}{2}]}\cdot s_i, r=1,2,\cdots ,j-1,$ where $[\alpha]$ denotes
the greatest integer equal to or smaller than $\alpha$.

$2)$ Otherwise, let $\varepsilon _k=s_k,$ i.e. no change is made
for other terms.

{\bf Lemma 2} \cite{Yang99,YHZ96,YZH96} \ \ Given a real
polynomial $f(x)=a_0x^n+a_1x^{n-1}+\cdots +a_n\in P^n,$ if the
number of sign changes in the revised sign list of its
discrimination sequence is $\nu$, and the number of nonzero
elements in the revised sign list of its discriminant sequence is
$\mu$, then the number of the distinct real roots of $f(x)$ is
$\mu-2\nu $.

{\bf Remark 1}\ \ The discrimination sequence of polynomial $f(x)$
can also be constructed by the principal sub-matrices of Bezout
matrix \cite{YHZ96,YZH96} of $f(x)$ and $f^{\prime }(x)$; the
number of distinct real roots of the polynomial $f(x)$ can also be
determined by the sign differences of Bezout matrix
\cite{YHZ96,YZH96} of $f(x)$ and $f^{\prime }(x)$ .

{\bf Remark 2}\ \ The original complete discrimination system of
polynomial \cite{Yang99,YHZ96,YZH96} is more general than Lemma 2,
which can also be used to determine the number of complex roots
and the multiplicities of repeated roots.

\vskip 0.7cm \section{The Global Hyperbolicity of Difference
Systems}

Consider the difference system

\begin{equation}
x(t)-\sum^{N}_{k=1}B_{k}x(t-\gamma_{k}\cdot r)=0 \label{eq31}
\end{equation}
and its characteristic equation

\begin{equation}
E(\lambda,r,B)=\det [I-\sum^{N}_{k=1}B_{k}e^{-\lambda \gamma_{k}\cdot r}]=0 \label{eq41}
\end{equation}
where all parameters are the same as before. In order to analyze
the global hyperbolicity of the difference system $(\ref{eq31})$,
for its characteristic equation $(\ref{eq41})$, by \cite{AH80}, we
have

{\bf Lemma 3} \cite{AH80} \ \ Suppose $r_1,r_2,\cdots,r_M$ are
commensurable, that is, there exist a $\beta >0$ and some integers
$n_k$ such that $r_k=n_k \beta,k=1,\cdots,M$, then there exists an
integer $p$, such that
 $E(\lambda,r,B)$ is a polynomial of some degree $p$ in  $e^{-\lambda
 \beta}$.
Denote the $p^{th}$-degree coefficient of $e^{-\lambda \beta}$ as
$A_p$, then

$$E(\lambda,r,B)=A_p \prod ^{p}_{\nu =1}(e^{- \lambda \beta}-s_{\nu})$$
and

$$\bar{Z}(B,r)=Z(B,r)=\{ -\displaystyle\frac{1}{\beta} \ln |s_{\nu} | ,\nu=1,2,\cdots,p \}$$
where $Z(B,r)=\{\mbox{Re}\lambda: E(\lambda,r,B)=0 \},
\bar{Z}(B,r)=\overline{Z(B,r)}= \mbox{cl}\{Z(B,r)\}$
is the closure of $Z(B,r)$.

Now, consider the corresponding difference system of System
$(\ref{eq6})$

\begin{equation}
y(t)-\sum^{N}_{k=1}B_{k}y(t-k \tau)=0 \label{eq610}
\end{equation}
Its characteristic equation is

\begin{equation}
e(\lambda,\tau,B)=\det [I-\sum^{N}_{k=1}B_{k}e^{-\lambda k
\tau}]=0 \label{eq620}
\end{equation}

{\bf Lemma 4} \ \ The difference system $(\ref{eq610})$ is
globally hyperbolic at $B$ if and only if
$$\forall \theta \in [0, 2\pi], \det[I-\sum^{N}_{k=1}B_{k}e^{i k \theta}]\neq 0.$$

{\bf Proof.} From Lemma 3, the number of the roots of the
characteristic equation $(\ref{eq620})$ is finite, and distribute
discretely in $R \times R$ plane. So the following conclusion
comes naturally:

If the real part of all characteristic roots is nonzero, then
there exists a $\delta >0 $, such that $\{\mbox{Re}\lambda:
E(\lambda,r,B)=0 \} \cap [-\delta,\delta]= \phi$.

That is, the difference system $(\ref{eq610})$ is globally
hyperbolic at $B$ if and only if
$$\forall y \in R, e(iy,\tau,B)=\det [I-\sum^{N}_{k=1}B_{k}e^{-iy k  \tau}]\neq 0
.$$

Letting $-y \tau = \theta,$ the lemma is proved. $\square$

Owing to the above lemma, it may be feasible to transform the
problem of global hyperbolicity of the difference system
$(\ref{eq610})$ into the real root existence of some related
polynomials.

\vskip 0.7cm
 \section{ Algebraic Criteria for
Delay-independent Stability of Neutral
       Differential Systems }

Consider the linear multi-delay neutral differential system
\begin{equation}
\dot{x}(t)-\sum^{N}_{k=1}B_{k}\dot{x}(t-k\tau)=A_{0}x(t)+\sum^{N}_{k=1}A_{k}x(t-k\tau) \label{eq66}
\end{equation}

The characteristic equation is
\begin{equation}
H(\lambda,\tau,A,B)= \det [\lambda
(I-\sum^{N}_{k=1}B_{k}e^{-\lambda k
\tau})-A_0-\sum^N_{k=1}A_ke^{-\lambda k \tau}]=0 \label{eq661}
\end{equation}

The corresponding difference system is
\begin{equation}
y(t)-\sum^{N}_{k=1}B_{k}y(t-k \tau)=0 \label{eq61}
\end{equation}
and its characteristic equation is

\begin{equation}
e(\lambda,\tau,B)=\det [I-\sum^{N}_{k=1}B_{k}e^{-\lambda k
\tau}]=0 \label{eq62}
\end{equation}

{\bf Theorem 1} \ \ The neutral differential system $(\ref{eq66})$
is delay-independent stable if and only if\\

(i)  $\det[I-\sum^{N}_{k=1}B_{k}e^{i k \theta}]\neq 0, \forall \theta \in [0,2\pi],$\\

(ii)  $\mbox{Re}\lambda[(I-\sum^{N}_{k=1}B_{k})^{-1}\sum^{N}_{k=0}A_{k}]<0,$\\

(iii)  $\det[iy(I-\sum^{N}_{k=1}B_{k}e^{i k
\theta})-A_{0}-\sum^{N}_{k=1}A_{k}e^{i k \theta}]\neq 0, \forall
\theta \in [0,2\pi],\forall y \in R,y \neq 0$.

{\bf Proof.} By lemma 1 and lemma 4, letting $-y \tau = \theta$,
the theorem is proved. $\square$

Condition (i) and condition (iii) are also transcendental, which
are still difficult to test numerically. Next, we will give a
simple criterion.

The bilinear transform $\omega =\displaystyle\frac {1+z}{1-z}$ is
a one-to-one mapping between the set $\{ \omega
=e^{i\theta}:\theta \in [0,2\pi]\}$ and the set $\{ z=iy:y\in
R\}$. By this transformation, condition (i) in Theorem 1 is
equivalent to
$$\det [I-\sum^{N}_{k=1}B_{k}(\displaystyle\frac{1+iz}{1-iz})^{k}]\neq 0$$
that is
$$\det[(1-iz)^{N}I-\sum^{N}_{k=1}B_{k}(1+iz)^{k}(1-iz)^{N-k}]\neq 0$$

Let
$$d(z)=\det [(1-iz)^{N}I-\sum^{N}_{k=1}B_{k}(1+iz)^{k}(1-iz)^{N-k}]$$
and $f(z)=\mbox{Re}[d(z)],g(z)=\mbox{Im}[d(z)],$ where $f(z),g(z)$
are real polynomials in $z$. Suppose
\begin{equation}
\left\{
\begin{array}{l}
f(z)=a_{0}z^{l_1}+a_{1}z^{l_1-1}+\cdots +a_{l_1}, \\
g(z)=b_{0}z^{m_1}+b_{1}z^{m_1-1}+\cdots +b_{m_1}\hbox{¡£}
\end{array}
\right. \label{eqaa}
\end{equation}
where $a_{0}\neq 0, b_{0}\neq 0$.

Then condition (i) in Theorem 1 is equivalent to the condition
that (\ref{eqaa}) has no common real roots.

By the same transformation, condition (iii) in Theorem 1 is
equivalent to
$$\det [iy(I-\sum^{N}_{k=1}B_{k}(\displaystyle\frac{1+iz}{1-iz})^{k})
-A_{0}-\sum^{N}_{k=1}A_{k}(\displaystyle\frac{1+iz}{1-iz})^{k}]\neq
0$$
namely,
$$\det [(iyI-A_{0})(1-iz)^{N}-iy\sum^{N}_{k=1}B_{k}(1+iz)^{k}(1-iz)^{N-k}
-\sum^{N}_{k=1}A_{k}(1+iz)^{k}(1-iz)^{N-k}]\neq 0$$

Let
$$D(z,y)=\det [(iyI-A_{0})(1-iz)^{N}-iy\sum^{N}_{k=1}B_{k}(1+iz)^{k}(1-iz)^{N-k}
-\sum^{N}_{k=1}A_{k}(1+iz)^{k}(1-iz)^{N-k}]$$
and
$F(z,y)=\mbox{Re}[D(z,y)],G(z,y)=\mbox{Im}[D(z,y)],$ where
$F(z,y),G(z,y)$ are real polynomials in $(z,y)$.
Suppose
\begin{equation}
\left\{
\begin{array}{l}
F(z,y)=a_{0}(y)z^{l_2}+a_{1}(y)z^{l_2-1}+\cdots +a_{l_2}(y), \\
G(z,y)=b_{0}(y)z^{m_2}+b_{1}(y)z^{m_2-1}+\cdots
+b_{m_2}(y)\hbox{.}
\end{array}
\right. \label{eqbb}
\end{equation}

Then condition (iii) in Theorem 1 is equivalent to the condition
that (\ref{eqbb}) has no common real roots $(z, y), z \in R, y\in
R-\{0\}$.

Before presenting the stability criteria, we first introduce two
lemmas in algebra theory \cite{Gan59,YZH96}.

{\bf Lemma 5}\ \ The equation $(\ref{eqaa})$ has at least one
common root if and only if its resultant equals zero, i.e.,
\begin{equation}
R(f,g)=\det \left[
\begin{array}{lllllll}
a_{0} & a_{1} & a_{2} & \cdots  & a_{l_1} &  &  \\
& a_{0} & a_{1} & \cdots  &  & a_{l_1} &  \\
&  & \cdots \cdots  &  &  &  &  \\
&  &  & a_{0} & a_{1} & \cdots  & a_{l_1} \\
b_{0} & b_{1} & b_{2} & \cdots  & b_{m_1} &  &  \\
& b_{0} & b_{1} & \cdots  &  & b_{m_1} &  \\
&  & \cdots \cdots  &  &  &  &  \\
&  &  & b_{0} & b_{1} & \cdots  & b_{m_1}
\end{array}
\right]
\begin{array}{l}
\left.
\begin{array}{l}
\\
\\
\\
\\
\end{array}
\right\} m_1 \mbox{ rows}  \\
\left.
\begin{array}{l}
\\
\\
\\
\\
\end{array}
\right\} l_1 \mbox{ rows}
\end{array}
 =0. \label{p1}
\end{equation}

{\bf Lemma 6}\ \ Suppose $a_{0}(y)\neq 0$ or $b_{0}(y)\neq 0$, the
equation $(\ref{eqbb})$ has at least one common root $(z,y)$ and
$y\in R$ if and only if its resultant
\begin{equation}
R(F,G)=\det \left[
\begin{array}{lllllll}
a_{0}(y) & a_{1}(y) & a_{2}(y) & \cdots  & a_{l_2}(y) &  &  \\
& a_{0}(y) & a_{1}(y) & \cdots  &  & a_{l_2}(y) &  \\
&  & \cdots \cdots  &  &  &  &  \\
&  &  & a_{0}(y) & a_{1}(y) & \cdots  & a_{l_2}(y) \\
b_{0}(y) & b_{1}(y) & b_{2}(y) & \cdots  & b_{m_2}(y) &  &  \\
& b_{0}(y) & b_{1}(y) & \cdots  &  & b_{m_2}(y) &  \\
&  & \cdots \cdots  &  &  &  &  \\
&  &  & b_{0}(y) & b_{1}(y) & \cdots  & b_{m_2}(y)
\end{array}
\right]
\begin{array}{l}
\left.
\begin{array}{l}
\\
\\
\\
\\
\end{array}
\right\} m_2 \mbox{ rows}  \\
\left.
\begin{array}{l}
\\
\\
\\
\\
\end{array}
\right\} l_2 \mbox{ rows}
\end{array}
=0 \label{p2}
\end{equation}
has a real root $y\in R$.

Similarly, $F(z,y),G(z,y)$ can be rewritten in the following
formula:
\begin{equation}
\left\{
\begin{array}{l}
F(z,y)=c_{0}(z)y^{l_3}+c_{1}(z)y^{l_3-1}+\cdots +c_{l_3}(z), \\
G(z,y)=d_{0}(z)y^{m_3}+d_{1}(z)y^{m_3-1}+\cdots
+d_{m_3}(z)\hbox{.}
\end{array}
\right. \label{eqcc}
\end{equation}
Then condition (iii) in Theorem 1 is equivalent to the condition
that (\ref{eqcc}) has no common real roots $(z, y), z \in R, y\in
R-\{0\}$.

Also from algebra theory \cite{Gan59,YZH96}, we have:

{\bf Lemma 7}\ \ Suppose $c_{0}(z)\neq 0$ or $d_{0}(z)\neq 0$, the
equation $(\ref{eqcc})$ has at least one common root $(z,y)$ and
$z\in R$ if and only if its resultant
\begin{equation}
\widetilde{R}(F,G)=\det \left[
\begin{array}{lllllll}
c_{0}(z) & c_{1}(z) & c_{2}(z) & \cdots  & c_{l_3}(z) &  &  \\
& c_{0}(z) & c_{1}(z) & \cdots  &  & c_{l_3}(z) &  \\
&  & \cdots \cdots  &  &  &  &  \\
&  &  & c_{0}(z) & c_{1}(z) & \cdots  & c_{l_3}(z) \\
d_{0}(z) & d_{1}(z) & d_{2}(z) & \cdots  & d_{m_3}(z) &  &  \\
& d_{0}(z) & d_{1}(z) & \cdots  &  & d_{m_3}(z) &  \\
&  & \cdots \cdots  &  &  &  &  \\
&  &  & d_{0}(z) & d_{1}(z) & \cdots  & d_{m_3}(z)
\end{array}
\right]
\begin{array}{l}
\left.
\begin{array}{l}
\\
\\
\\
\\
\end{array}
\right\} m_3 \mbox{ rows}  \\
\left.
\begin{array}{l}
\\
\\
\\
\\
\end{array}
\right\} l_3 \mbox{ rows}
\end{array}
=0 \label{p3}
\end{equation}
 has a real root $z\in R$.

From Lemma 6 and Lemma 7, we have:

{\bf Lemma 8}\ \ If $(z_{0},y_{0})$ are a pair of real roots of
$F(z,y)=0, G(z,y)=0$, then $R(F,G)(y_{0})=0,
\widetilde{R}(F,G)(z_{0})=0$.

{\bf Proof.} The lemma is proved by using lemma 6 and lemma 7
directly. $\square$

{\bf Remark 4}\ \ It should be noticed that the converse
proposition of Lemma 8 is not correct, which can be tested by some
examples easily.

Now we are in position to present our main result. The algebraic
criteria for the delay-independent stability of the neutral
differential system $(\ref{eq66})$ are as follows:

{\bf Theorem 2} \ \ The neutral differential system $(\ref{eq66})$
is delay-independent stable if and only if

(i) $(\ref{eqaa})$ has no common real roots.

(ii) $\mbox{Re}\lambda[(I-\sum^{N}_{k=1}B_{k})^{-1}\sum^{N}_{k=0}A_{k}]<0.$

(iii) $(\ref{eqbb})$(or $(\ref{eqcc})$) has no common real roots
$(z,y),z\in R,y\in R - \{ 0 \}.$

{\bf Proof.} From the above discussion, condition (i) in theorem 1
is equivalent to $(\ref{eqaa})$ having no common real roots;
condition (ii) in theorem 2 is same as condition (ii) in theorem
1; condition (iii) in theorem 1 is equivalent to $(\ref{eqbb})$
(or $(\ref{eqcc})$) having no common real roots $(z,y),z\in R,y\in
R - \{ 0 \}$. This completes the proof. $\square$

By Lemma 5, the condition (i) above \textit{is equivalent to} one
of the following conditions:

(a1) the determinant $R(f,g) \neq 0.$

(a2) $R(f,g)=0$ and $g(z)\neq 0$, where $z$ is the real root of
$f(z)$, and $R(f,g)$ is the resultant of $(\ref{eqaa}).$

(a3) $R(f,g)=0$ and $f(z)\neq 0$, where $z$ is the real root of
$g(z)$, and $R(f,g)$ is the resultant of $(\ref{eqaa}).$

Similarly, by Lemma 6, the condition (iii) above \textit{is
equivalent to} one of the following conditions:

(b1) when $a_{0}(y)\neq 0$ or $b_{0}(y)\neq 0,$
 $R(f,g)$ has no nonzero real root, where $R(f,g)$ is the resultant of $(\ref{eqbb}).$

(b2) when $a_{0}(y)\neq 0$ or $b_{0}(y)\neq 0,$
 $R(f,g)$ has nonzero real root. But when taking the root back into $(\ref{eqbb}),$
  then $(\ref{eqbb})$ has no common real root.

(b3) when there exists a real number $y$, such that $a_{0}(y)= 0$
and $b_{0}(y)= 0,$ then take this real number $y$ into
$(\ref{eqbb}),$ and $(\ref{eqbb})$ has no common real root.

Dually, by Lemma 7, the condition (iii) above \textit{is also
equivalent to} one of the following conditions:

(c1) when $c_{0}(z)\neq 0$ or $d_{0}(z)\neq 0,$
 $\widetilde{R}(f,g)$ has no nonzero real root, where $\widetilde{R}(f,g)$ is the resultant of $(\ref{eqcc}).$

(c2) when $c_{0}(z)\neq 0$ or $d_{0}(z)\neq 0,$
 $\widetilde{R}(f,g)$ has nonzero real root. But when taking the root back into $(\ref{eqcc}),$
  then $(\ref{eqcc})$ has no common real root.

(c3) when there exists a real number $z$, such that $c_{0}(z)= 0$
and $d_{0}(z)= 0,$ then take this real number $z$ into
$(\ref{eqcc}),$ and $(\ref{eqcc})$ has no common real root.

Furthermore, by Lemma 8, the condition (iii) above can also be
checked with the following condition:

(d1) $\forall y_{0} \in \{y|R(F,G)(y)=0\}\cap\{R-\{0\}\}, \forall
z_{0} \in \{z|\widetilde{R}(F,G)(z)=0\}\cap R, (y_{0}, z_{0})$ are
not the roots of Equations $(\ref{eqbb})$ or Equations
$(\ref{eqcc}).$

{\bf Remark 5} \ \ The conditions in Theorem 2 are all algebraic
conditions. Condition (ii) can be checked by the well-known
Hurwitz Criterion \cite{Gan59}. Determining real roots of
polynomials in conditions (i) and (iii) in Theorem 2 can be
carried out by the ``Complete Discrimination System for
Polynomials" mentioned before.

{\bf Remark 6} \ \ The conditions in Theorem 2 are necessary and
sufficient. The sign list of the discrimination sequence of
polynomials with symbolic coefficients can be obtained easily by
computer \cite{Yang99,YHZ96,YZH96}. Therefore, ``on-line"
determining the delay-independent stability of neutral
differential systems $(\ref{eq66})$ can be realized, namely, an
efficient algorithm can be set up by this theorem.

{\bf Remark 7} \ \ System $(\ref{eq66})$ is very general, which
covers various forms of systems studied in
\cite{CL95,GTY01,Han01,HJZ84,HH96,Kam82,Li88,PW00,RL02,RZ99,VN97,WLB02,Wu94,YW00}.
Theorem 1 and Theorem 2 generalize the relevant results in
\cite{GTY01,YW00}. More specifically, if $B_k=0,k=1,\cdots,N,$ in
$(\ref{eq66})$, then System $(\ref{eq66})$ degenerates to a
retarded differential system. The results in this paper  are
completely consistent with that in \cite{GTY01,YW00}.

Just like the discussion in \cite{GTY01}, when System
$(\ref{eq66})$ is not delay-independent stable, the maximal delay
bound guaranteeing stability can be determined as following
\begin{equation}
T:=\min \left\{\tau |\tau >0,(z,y)\in \left\{
\begin{array}{l}
F(z,y)=0,\\
G(z,y)=0,
\end{array},
z\in R,y\in R-\{0\} \right\}, e^{i\theta
}=\frac{1+iz}{1-iz},-y\tau =\theta \right \}. \label{eqTT}
\end{equation}

In the following section, some numerical examples are provided to
to test the delay-independent stability of delay systems and
compute the maximal delay bound when the systems are not
delay-independent stable.

\vskip 0.7cm
\section{ Corollaries and Examples}

In what follows, we will present some explicit algebraic criteria
for some simple neutral differential systems.

{\bf Corollary 1} \ \ The system
\begin{equation}
\dot{x}(t)+c\dot{x}(t-\tau )+ax(t)+bx(t-\tau )=0 \quad  (1+c \neq
0) \label{eqrz99}
\end{equation}
is delay-independent stable if and only if
$$(1-c)(b-a)<0 \hbox{ and } (1+c)(b+a)>0$$

{\bf Proof.} Let

$f(z)=\mbox{Re}[(1-iz)+c(1+iz)]=\allowbreak c+1$

$g(z)=\mbox{Im}[(1-iz)+c(1+iz)]=\allowbreak cz-z$

$[\lambda +(1+c)^{-1}(a+b)]=0, \lambda =-(1+c)^{-1}(a+b)$

$F(z,y)=\mbox{Re}[(iy+a)(1-iz)+iyc(1+iz)+b(1+iz)]=\allowbreak
yz\left( 1-c\right) +a+b$

$G(z,y)=\mbox{Im}[(iy+a)(1-iz)+iyc(1+iz)+b(1+iz)]=\allowbreak
y\left( c+1\right) +z\left( b-a\right) $

Obviously, $f(z)=0$ and $g(z)=0$ have no common real root.

By the condition (ii) in theorem 2, we have $\lambda
=-(1+c)^{-1}(a+b)<0,$ that is $(1+c)(b+a)>0.$

Finally, it is easy to see that $F(z,y)=0$ and $G(z,y)=0$ has no
common real roots $(z,y),z\in R,y\in R-\{0\},$ if and only if
$(1-c)(b-a)<0.$ Thus by theorem 2, we complete the proof.

{\bf Corollary 2} \ \ The system
\begin{equation}
\dot{x}(t)+c\dot{x}(t-2\tau )+ax(t)+ax(t-\tau )=0 \quad (1+c \neq
0) \label{eqrl02}
\end{equation}
is delay-independent stable if and only if
$$a>0 \hbox{ and }-1< c \leq \frac{1}{3}$$

{\bf Proof.} Let

$f(z)=\mbox{Re}[(1-iz)^{2}+c(1+iz)^{2}]=\allowbreak c\left(
1-z^{2}\right) -z^{2}+1$

$g(z)=\mbox{Im}[(1-iz)^{2}+c(1+iz)^{2}]=\allowbreak 2cz-2z$

$[\lambda -(1+c)^{-1}(-a+-a)] = 0, \lambda =-2a(1+c)^{-1}$

$F(z,y)=\mbox{Re}[(iy+a)(1-iz)^{2}+iyc(1+iz)^{2}+a(1+iz)(1-iz)]=\allowbreak
a+2yz-2cyz+az^{2}+a\left( 1-z^{2}\right)$

$G(z,y)=\mbox{Im}[(iy+a)(1-iz)^{2}+iyc(1+iz)^{2}+a(1+iz)(1-iz)]=\allowbreak
y\left( 1-z^{2}\right) -2az+cy\left( 1-z^{2}\right) $

Obviously, $f(z)=0$ and $g(z)=0$ have no common real root.

By the condition (ii) in theorem 2, we have $\lambda
=-2a(1+c)^{-1}<0$ that is $a>0,1+c>0;$ or $a<0,1+c<0.$

Finally, the resultant of $F(z,y)$ and $G(z,y)$

$\widetilde{R}(F,G)(z)\det\left[
\begin{array}{cc}
2z-2cz & 2a \\
(1-z^{2})(1+c) & -2az%
\end{array}%
\right] = \allowbreak 2a(3c-1)z^{2}-2a(1+c)$

It is easy to see that $F(z,y)=0$ and $G(z,y)=0$ has no common
real roots $(z,y),z\in R,y\in R-\{0\},$ if and only if $a>0$ and
$-1<c \leq \displaystyle\frac{1}{3}.$ Thus by theorem 2, we
complete the proof.

{\bf Remark 8} \ \ Corollary 1 and Corollary 2 are consistent with
the relevant results in \cite{RZ99,RL02}, respectively. They have
also discussed delay-dependent stability of System
$(\ref{eqrz99})$ and System $(\ref{eqrl02})$.

Next, we will present some numerical examples using the main
results of this paper.

 {\bf Example 1} \ \ Consider the following system
$$\dot{X}(t)-C\dot{X}(t-\tau )=AX(t)+BX(t-\tau )$$
where
$$C=\left[
\begin{array}{ll}
0.1 & 0 \\
0 & 0.1%
\end{array}
\right] ,A=\left[
\begin{array}{ll}
-3 & -2 \\
1 & 0%
\end{array}
\right] ,B=\alpha \left[
\begin{array}{ll}
0 & 1 \\
1 & 0%
\end{array}
\right]  $$
and $\alpha$ is a nonzero constant.

We now consider the maximal stability bound in terms of $\alpha$.
Let $I$ be the $2 \times 2$ identity matrix.

$f(z):=\mbox{Re}[\det [I(1-iz)-C(1+iz)]]=\allowbreak
0.\,81-1.\,21z^2$

$g(z):=\mbox{Im}[\det [I(1-iz)-C(1+iz)]]=\allowbreak -1.\,98z$

$\det [\lambda I-(I-C)^{-1}(A+B)]=\allowbreak \lambda
^{2}+3.\,3333\lambda +2.\,4691+\allowbreak 1.\,2346\alpha
-1.\,2346\alpha ^{2}$

$F(z,y):=\mbox{Re}[\det
[(iyI-A)(1-iz)-iyC(1+iz)-B(1+iz)]]=\allowbreak y^{2}\left(
1.\,\allowbreak 21z^{2}-0.81\right) +6.0yz+\alpha -2\allowbreak
z^{2}-\alpha ^{2} +z^{2}\alpha +z^{2}\alpha ^{2}+2\allowbreak
\allowbreak $

$G(z,y):=\mbox{Im}[\det
[(iyI-A)(1-iz)-iyC(1+iz)-B(1+iz)]]=\allowbreak 1.\,\allowbreak
98y^{2}z+y\left( 2.\,\allowbreak 7-3.\,\allowbreak 3z^{2}\right)
-4z-\allowbreak 2z\alpha ^{2}\allowbreak $

Obviously, $f(z)$ and $g(z)$ have no common real roots.

If $\mbox{Re}\lambda[(I-C)^{-1}(A+B)]<0,$ by the Hurwitz
criterion, we have $-0.\,99997<\alpha<2.0$.

Finally, the resultant of $F(z,y)$ and $G(z,y)$

$\widetilde{R}(F,G)(z) :=( -5.\,9049\alpha +5.\,9049\alpha
^2-11.\,81) +( -57.\,802-24.\,592\alpha ^2+12.\,96\alpha ^4
-14.\,256\alpha ^3-11.\,875\alpha ) z^2+( 6.\,336\alpha
-31.\,68\alpha ^4-81.\,081\alpha ^2-106.\,01 +3.\,168\alpha ^3)
z^4+( 25.\,483\alpha -40.\,608\alpha ^2-86.\,346+17.\,424\alpha ^3
+19.\,36\alpha ^4) z^6 +( 13.\,177\alpha +13.\,177\alpha
^2-26.\,354) z^8$

By a careful calculation, $\widetilde{R}(F,G)(z)=0$ has no real
root when $\alpha \in (-0.99997, 1]$. In fact, if $\alpha \in
(-0.99997, 1]$, all coefficients of $\widetilde{R}(F,G)(z)$ have
the same sign. Thus, $F(z,y)=0$ and $G(z,y)=0$ have no common real
roots. Hence, by Theorem 2, when $\alpha \in (-0.99997, 1]$, the
system is delay-independent stable.

{\bf Remark 9} \ \ The criteria in \cite{HH96,Li88} do not work
here. \cite{PW00} concluded that the system is delay-independent
stable when $|\alpha |\leq 0.989$. By using the criteria in this
paper, we conclude that the system is delay-independent stable
when $\alpha \in (-0.99997, 1]$, which gives a less conservative
bound of $\alpha$.

{\bf Example 2}\ \ Consider the following system
$$\dot{X}(t)-C\dot{X}(t-\tau )=AX(t)+BX(t-\tau )$$
where
$$C=\left[
\begin{array}{ll}
0 & 0.4 \\
0.4 & 0%
\end{array}
\right] ,A=\left[
\begin{array}{ll}
-1 & 0 \\
0 & -1%
\end{array}
\right] ,B=\alpha \left[
\begin{array}{ll}
0 & 1 \\
1 & 0%
\end{array}
\right] $$
and $\alpha$ is a nonzero constant.

We now consider the maximal stability bound in terms of $\alpha$.
Let $I$ be the $2 \times 2$ identity matrix.

$f(z):=\mbox{Re}[\det [I(1-iz)-C(1+iz)]]=\allowbreak
0.\,84-0.\,84z^2$

$g(z):=\mbox{Im}[\det [I(1-iz)-C(1+iz)]]=\allowbreak -2.\,32z$

$\det [\lambda I-(I-C)^{-1}(A+B)]]=\allowbreak \lambda
^2+2.\,381\lambda -0.\,95238\lambda \alpha +\allowbreak
1.\,1905-1.\,1905\alpha ^2$

$F(z,y):=\mbox{Re}[\det
[(iyI-A)(1-iz)-iyC(1+iz)-B(1+iz)]]=\allowbreak
-0.\,84y^2+1+4zy+0.\,84z^2y^2-z^2 +1.\,6y\alpha z-\allowbreak
\alpha ^2+\alpha ^2z^2$

 $G(z,y):=\mbox{Im}[\det
[(iyI-A)(1-iz)-iyC(1+iz)-B(1+iz)]]=\allowbreak
2y+2.\,32zy^2-2z-2z^2y-0.\,8y\alpha +0.\,8z^2y\alpha -\allowbreak
2\alpha ^2z$

Obviously, $f(z)$ and $g(z)$ have no common real roots.

If $\mbox{Re}\lambda[(I-C)^{-1}(A+B)]<0,$ by the Hurwitz
criterion, we have $-1.0<\alpha <1.0.$

Finally, the resultant of $F(z,y)$ and $G(z,y)$

$R(F,G)(y):=\allowbreak ( -4.0+8.0\alpha ^4-4.0\alpha ^8)+ (
2.\,56\alpha ^2-2.\,56\alpha ^6-16.0+16.0\alpha ^4) y^2 +(
7.\,3856\alpha ^4+5.\,12\alpha ^2-23.\,795) y^4+( 2.\,4945\alpha
^2-15.\,59) y^6-3.\,7978y^8$

It is easy to see that $R(F, G)(y)=0$ has no real root when
$\alpha \in (-1.0, 1.0)$. In fact, if $\alpha \in (-1.0, 1.0)$,
all coefficients of $R(F,G)(z)$ have the same sign. Thus,
$F(z,y)=0, G(z,y)=0$ have no common real roots. Hence, when
$\alpha \in (-1.0, 1.0)$, the system is delay-independent stable.

{\bf Remark 10} \ \ The criteria in \cite{HH96,Li88,PW00} can work
here, and they can get $\alpha $'s bounds respectively as: $
|\alpha | \leq 0.2$ \cite{HH96}; $ |\alpha | \leq 0.2$
\cite{Li88};  $ |\alpha | \leq 0.9165$ \cite{PW00}. The result of
this paper is $\alpha \in (-1.0, 1.0)$, which obviously gives a
less conservative bound of $\alpha$.

{\bf Example 3}\ \ Consider the following system
$$\dot{X}(t)-C\dot{X}(t-\tau )=AX(t)+BX(t-\tau )$$
where\\
$$C=\left[
\begin{array}{ll}
0.2 & 0 \\
0 & 0.2%
\end{array}
\right] ,A=\left[
\begin{array}{ll}
-2 & 0 \\
0 & -1%
\end{array}
\right] ,B=\left[
\begin{array}{ll}
0 & 0.5 \\
0.5 & 0%
\end{array}
\right] $$

We will discuss whether the system is delay-independent stable.
Let $I$ be the $2 \times 2$ identity matrix.

$f(z):=\mbox{Re}[\det [I(1-iz)-C(1+iz)]]=\allowbreak
0.\,64-1.\,44z^2$

$g(z):=\mbox{Im}[\det [I(1-iz)-C(1+iz)]]=\allowbreak -1.\,92z$

$\det [\lambda I-(I-C)^{-1}(A+B)]]=\allowbreak \lambda
^{2}+3.\,75\lambda +2.\,7344=0$

$F(z,y):=\mbox{Re}[\det
[(iyI-A)(1-iz)-iyC(1+iz)-B(1+iz)]]=\allowbreak
-0.\,64y^2+6.0zy+1.\,75+1.\,44z^2y^2-\allowbreak 1.\,75z^2$

$G(z,y):=\mbox{Im}[\det
[(iyI-A)(1-iz)-iyC(1+iz)-B(1+iz)]]=\allowbreak
2.\,4y+1.\,92zy^2-4.\,5z-3.\,6z^2y$

Obviously, $f(z)$ and $g(z)$ have no common real roots.

By the Hurwitz criterion,  it is obvious $Re \lambda
[I-(I-C)^{-1}(A+B)] <0.$

The resultant of $F(z,y)$ and $G(z,y)$, $R(F,G)(y),$ is

$R(F,G)(y):=-35.\,684y^6-122.\,57y^4-\allowbreak
3.\,3974y^8-152.\,46y^2-\allowbreak 62.\,016$

$R(F,G)(y)=0$ has no real roots, that is, $F(z,y)=0,G(z,y)=0$ have
no common real roots.

Hence, the system is delay-independent stable.

{\bf Remark 11} \ \ By the criteria in \cite{KY88,PW00}, the
maximum allowable bound guaranteeing asymptotic stability of the
system was obtained, they are $\tau_{max}=0.1352$ in \cite{KY88}
and $\tau_{max}=0.7516$ in \cite{PW00}, respectively. However,
this paper concludes that the system is delay-independent stable.
The reason for the different conclusions is that the criterion in
this paper is necessary and sufficient, whereas the criteria in
\cite{KY88,PW00} are only sufficient conditions.

The above three examples were all discussed in \cite{PW00}.
Sharper stability bounds are obtained by using the algebraic
criteria in Theorem 2.

 {\bf Example 4}\ \ Consider the system
$$\dot{X}(t)-C\dot{X}(t-\tau)=AX(t)+BX(t-\tau)$$
where
$$C=\left[
\begin{array}{lll}
0.5 & 0 & 0 \\
0 & 0.2 & 0 \\
0 & 0 & 0.3
\end{array}
\right] ,A=\left[
\begin{array}{lll}
-3 & -2 & -2 \\
2 & -2 & -2 \\
0 & 0 & -2
\end{array}
\right] ,B=\left[
\begin{array}{lll}
-1 & 1 & 1 \\
1 & -1 & -1 \\
0 & 0 & -1
\end{array}
\right]$$

We will discuss whether the system is delay-independent stable.
Let $I$ be the $3 \times 3$ identity matrix.

$f(z):=\mbox{Re}[\det[I(1-iz)-C(1+iz)]]=\allowbreak
0.\,28-3.\,6z^2$

$g(z):=\mbox{Im}[\det[I(1-iz)-C(1+iz)]]=\allowbreak
-1.\,78z+2.\,34z^3$

$\det [\lambda I-(I-C)^{-1}(A+B)]=\allowbreak \lambda
^3+16.\,036\lambda ^2+87.\,857\lambda +\allowbreak 160.\,71$

$F(z,y):=\mbox{Re}[\det
[(iyI-A)(1-iz)-iyC(1+iz)-B(1+iz)]]=\allowbreak 45+72.\,4zy
+24.\,75z^2y^2-1.\,78y^3z+\allowbreak
2.\,34z^3y^3-10.\,4z^3y-4.\,49y^2-\allowbreak 35z^2$

$G(z,y):=\mbox{Im}[\det
[(iyI-A)(1-iz)-iyC(1+iz)-B(1+iz)]]=\allowbreak 19.\,89y^2z
-52.\,6z^2y+3.\,6y^3z^2-\allowbreak
6.\,87z^3y^2-0.\,28y^3+5z^3-75z+24.\,6y$

It is easy to see that $f(z)=0,g(z)=0$ have no common real roots.

By the Hurwitz criterion, all roots of $\det [\lambda
I-(I-C)^{-1}(A+B)]$ have negative real parts.

Furthermore, by a careful calculation using Lemma 8, $F(z,y)$ and
$G(z,y)$ have no common real roots.

Hence, the system is delay-independent stable.

{\bf Remark 12} Let $C=0$ in this example, then it becomes a
retarded differential system. In this case, the conclusion in this
paper is still correct and is consistent with the one in
\cite{GTY01,YW00}.

Finally, we present a more complex example on computing the
maximal allowable delay bound.

{\bf Example 5}\ \ Consider the following system
$$\dot{X}(t)-C\dot{X}(t-\tau )=AX(t)+B_1X(t-\tau )+B_2X(t-2\tau )+B_3X(t-3\tau )$$
where
\begin{center}
$C=\left[
\begin{array}{llll}
0.02 & 0 & 0.03 & 0 \\
0 & 0.01 & 0 & 0 \\
0 & 0 & 0.5 & 0 \\
0 & 0 & 0 & 0.25
\end{array}
\right], A=\left[
\begin{array}{llll}
0 & 1 & 0 & 0 \\
0 & 0 & 1 & 0 \\
0 & 0 & 0 & 1 \\
-2 & -3 & -5 & -2
\end{array}
\right],$

$B_1=\left[
\begin{array}{llll}
-0.05 & 0.005 & 0.25 & 0 \\
0.005 & 0.005 & 0 & 0 \\
0 & 0 & 0 & 0 \\
-1 & 0 & -0.5 & 0
\end{array}
\right], B_2=\left[
\begin{array}{llll}
0.005 & 0.0025 & 0 & 0 \\
0 & 0 & 0.05 & 0 \\
0 & 0 & 0 & 0.0005 \\
-1 & -0.5 & -0.5 & 0
\end{array}
\right],$

$B_3=\left[
\begin{array}{llll}
0.0375 & 0 & 0.075 & 0.125 \\
0 & 0.05 & 0.05 & 0 \\
0.05 & 0.05 & 0 & 0 \\
0 & -2.5 & 0 & -1
\end{array}
\right]$
\end{center}

Let $I$ be the $4 \times 4$ identity matrix.

$f(z):=\mbox{Re}[\det [I(1-iz)-C(1+iz)]]=\allowbreak
0.\,36383-5.\,7048z^2+1.\,9316z^4$

$g(z):=\mbox{Im}[\det [I(1-iz)-C(1+iz)]]=\allowbreak
-2.\,4477z+5.\,5521z^3$

$\det [\lambda I-(I-C)^{-1}(A+B_1+B_2+B_3)]=\allowbreak \lambda
^4+3.\,949\lambda ^3 +16.\,661\lambda ^2+\allowbreak 19.\,9\lambda
+11.\,646$

$F(z,y):=\mbox{Re}[\det
[(iyI-A)(1-iz)^3-iyC(1+iz)(1-iz)^2-B_1(1+iz)(1-iz)^2
-B_2(1+iz)^2(1-iz)-B_3(1+iz)^3]]=\allowbreak
365.\,27y^2z^2+4587.\,5y^2z^6 -\allowbreak
771.\,21z^3y-2331.\,4z^7y+\allowbreak 1.\,9316z^{12}y^4
-911.\,59y^4z^6+\allowbreak 625.\,81y^4z^8-4.\,821z^{12}y^2
-\allowbreak 35.\,473y^4z^2+368.\,62y^4z^4+\allowbreak
26.\,68z^{11}y^3 -104.\,21z^{10}y^4+\allowbreak
2.\,008z^{12}-21.\,357z^{11}y +\allowbreak
321.\,27z^{10}y^2+1189.\,5z^4-\allowbreak 2558.\,4y^2z^4
+64.\,845zy+\allowbreak 2334.\,2z^5y-6.\,0618y^2-\allowbreak
211.\,53z^2 -1765.\,1z^6-\allowbreak
2396.\,7z^8y^2-555.\,07z^9y^3+\allowbreak 621.\,06z^9y
-17.\,459y^3z+\allowbreak
335.\,39y^3z^3-1413.\,6y^3z^5+\allowbreak 1747.\,8z^7y^3
+813.\,91z^8-\allowbreak 111.\,95z^{10}+.\,36383y^4+4.\,2371$

$G(z,y):=\mbox{Im}[\det
[(iyI-A)(1-iz)^3-iyC(1+iz)(1-iz)^2-B_1(1+iz)(1-iz)^2
-B_2(1+iz)^2(1-iz)-B_3(1+iz)^3]]=\allowbreak
-5.\,3583y^4z+140.\,1y^4z^3 +\allowbreak
1835.\,7y^3z^6-311.\,72y^4z^9+\allowbreak 97.\,231y^3z^2
+4002.\,2y^2z^5-\allowbreak
805.\,43y^3z^4+1542.\,3z^4y-\allowbreak 59.\,037z^{11}y^2
-278.\,51z^2y-\allowbreak 163.\,01z^{10}y+7.\,2403y-\allowbreak
44.\,703z -682.\,09y^4z^5-\allowbreak
1171.\,7y^2z^3+1457.\,7z^8y-\allowbreak 1190.\,1z^8y^3
-3875.\,1z^7y^2+\allowbreak
1062.\,4z^9y^2-2694.\,5z^6y+\allowbreak 69.\,443y^2z
+0.\,66541z^{12}y+\allowbreak
887.\,73z^7y^4-1.\,6155z^{12}y^3+\allowbreak 21.\,005z^{11}y^4
+164.\,45z^{10}y^3+\allowbreak
22.\,407z^{11}-1.\,4368y^3+\allowbreak 605.\,71z^3-1683.\,4z^5
+\allowbreak 1381.\,2z^7-356.\,38z^9$

Obviously, $f(z)=0,g(z)=0$ have no common real roots.

By Hurwitz Criterion, $\det [\lambda I-(I-C)^{-1}(A+B_1+B_2+B_3)]$
has only negative-real-part roots.

Further, by Lemma 8, we have the common real roots of $F(z,y)=0$
and $G(z,y)=0$ are $\{ y=-0.\,86798,z=10.\,823\} , \allowbreak \{
y=0.\,58092,z=10.\,7\} , \allowbreak \{ z=-1.\,6587,y=1.\,3876\} ,
\allowbreak \{ y=1.\,1338,z=-0.\,42313\} , \allowbreak \{
y=-0.\,58092,$ $z=-10.\,7\} , \allowbreak \{
y=0.\,31594,z=-1.\,2099\} , \allowbreak \{
z=-10.\,823,y=0.\,86798\} , \allowbreak \{
y=-1.\,1338,z=0.\,42313\} , \allowbreak \{
y=-3.\,3266,z=0.\,24369\} , \allowbreak \{
y=3.\,3266,z=-0.\,24369\} , \allowbreak \{
y=-1.\,3876,z\allowbreak =1.\,6587\} , \allowbreak \{
y=-0.\,31594,z =1.\,2099\}.$

Therefore, the system is not delay-independent stable by Theorem
2. The maximal delay bound guaranteeing the stability is
$T=0.14371$ (by taking all common real roots of $F(z,y)=0$ and
$G(z,y)=0$ into $(\ref{eqTT})$).

{\bf Remark 13} \ \ Set $C=0$ in Example 5, then the system
degenerates to a retard differential system. The maximal delay
bound was obtained in \cite{GTY01} as $T=0.4777$. The same
conclusion can be obtained by using the results in this paper.

From the above corollaries and examples, it is obvious that our
criteria work well on judging the delay-independent stability of
neutral differential systems, and its simplicity, accuracy,
convenience, and wide applicability greatly facilitate the
engineering practice.

\section{Conclusion}

This paper establishes some algebraic criteria for determining the
delay-independent stability of a class of neutral differential
systems
$\dot{x}(t)-\sum^{N}_{k=1}B_{k}\dot{x}(t-k\tau)=A_{0}x(t)+\sum^{N}_{k=1}A_{k}x(t-k\tau)$,
and presents a method for determining the maximal delay bound
guaranteeing stability if the systems are not delay-independent
stable. To the best of our knowledge, the algebraic criteria has
some noteworthy characters comparing with the related literature.

$\bullet$ {\bf Generality}

The system
$\dot{x}(t)-\sum^{N}_{k=1}B_{k}\dot{x}(t-k\tau)=A_{0}x(t)+\sum^{N}_{k=1}A_{k}x(t-k\tau)$
is more general, and covers various forms of systems studied in
\cite{CL95,GTY01,Han01,HJZ84,HH96,Kam82,Li88,PW00,RL02,RZ99,VN97,WLB02,Wu94,YW00}.
The criteria in this paper can be applied to all these systems,
and the relevant results  are consistent with
 or better than that of those literature, which have been illustrated by corollaries and examples.

$\bullet$ {\bf Necessary-Sufficient Conditions}

The criteria for delay-independent stability are necessary and
sufficient, and the broader conditions can be got comparing with
that of some literature (see example 1-3).

$\bullet$ {\bf Practicality}

The conditions in the criteria are all algebraic conditions.
Condition (ii) can be checked by Hurwitz Criterion \cite{Gan59}.
Determining real roots of polynomials in conditions (i) and (iii)
can be carried out by the ``Complete Discrimination System for
Polynomials" mentioned before.The sign list of the discrimination
sequence of polynomials with symbolic coefficients can be obtained
easily by computer \cite{Yang99,YHZ96,YZH96}. Therefore,
``on-line" operation can be realized.

\vskip 20pt
\vspace*{1\baselineskip}


\begin{thebibliography}{}

\bibitem{AH80}  {Avellar, C. E. and Hale, J. K., On the zeros of exponential polynomials,
\textit{J. Math. Anal. Appl.}, 1980, \textbf{73:} 434-452.}

\bibitem{BC63} {Bellman, R. and Cooke, K. L., {\it Differential-Difference Equations,}
Academic Press, New York, 1963.}

\bibitem{BL02} {Boukas, E. K. and Liu, Z. K., {\it Deterministic and Stochastic Time Delay System,}
Birkh$\ddot{a}$user, Boston, Basel, Berlin, 2002.}

\bibitem{Bru70}{Brumley, W. E., On the asymptotic behavior of solutions of differential-difference
equations of neutral type, \textit{J. Diff. Eqs.}, 1970,
\textbf{7:} 175-188.}

\bibitem{CL95}{Chen, J. and Latchman, H. A., Frequency sweeping tests for stability independent of delay,
\textit{IEEE Trans. Automat. Contr.}, 1995, \textbf{40:}
1640-1645.}

\bibitem{Gan59} {Gantmacher, F., {\it The Theory of Matrices,} Chelsea, New York, 1959.}

\bibitem{GTY01}  {Gu, N., Tan, M. and Yu, W. S., An algebra
test for unconditional stability of linear delay systems, \textit{
Proc. of the 40th IEEE Conference on Decision and Control (CDC
2001),} Orlando, Florida, USA, 2001, 4746-4747.}

\bibitem{Hal77} {Hale, J. K., {\it Theory of Functional Differential Equations,}
Springer-Verlag, New York, 1977.}

\bibitem{HIT85} {Hale, J. K., Infante, E. F. and Tsen, F. -S. P., Stability in linear delay equations,
\textit{J. Math. Anal. Appl.}, 1985, \textbf{105:} 533-555.}

\bibitem{HV93} {Hale, J. K. and Verduyn Lunel, S. M., {\it Introduction to Functional Differential Equations,}
Springer-Verlag, New York, 1993.}

\bibitem{Han01} {Han, Q. -L., On delay-dependent stability for neutral
delay-differential systems, \textit{ Int. J. Appl. Math. Comput. Sci.},
2001, \textbf{11:} 965-976.}

\bibitem{HJZ84} {Hertz, D., Jury, E. J. and Zeheb, E., Stability independent and
dependent of delay for delay differential systems, \textit{ J.
Franklin Inst.}, 1984, \textbf{318:} 143-150.}

\bibitem{HH96} {Hu, G. Di and Hu, G. Da, Some simple stability criteria of neutral delay-differential systems,
\textit{Appl. Math. Comput.}, 1996, \textbf{80:} 257-271.}

\bibitem{Kam82} {Kamen, E. W., Lineal systems with
commensurate time-delays: Stability and stabilization independent
of delay, \textit{IEEE Trans. Automat. Contr.}, 1982, \textbf{27:}
367-375.}

\bibitem{Kha99} {Kharitonov, V. L., Robust stability analysis of time delay systems: A survey,
\textit{Annual Reviews in Control}, 1999, \textbf{23:} 185-196.}

\bibitem{KTO03} {Kharitonov, V. L., Torres-Mu\~{n}oz, J. A. and Ortiz-Moctezuma, M. B., Polytopic families of
quasi-polynomials: vertex-type stability conditions, \textit{IEEE Trans. Circuits Syst., Part I}, 2003,
\textbf{CAS-50:} 1413-1420.}

\bibitem{KZ94} {Kharitonov, V. L. and Zhabko, A. P., Robust stability of time-delay systems,
\textit{IEEE Trans. Automat. Contr.}, 1994, \textbf{39:}
2388-2397.}

\bibitem{KY88} {Khusainov, D. Ya. and Yun'kova, E. V., Investigation of the stability of linear systems of
neutral type by the Lyapunov function method, \textit{Diff.
Uravn.}, 1988, \textbf{24:} 613-621.}

\bibitem{Li88} {Li, M. L., Stability of linear neutral delay-differential systems, \textit{
Bull. Austral. Math. Soc.}, 1988, \textbf{38:} 339-344.}


\bibitem{PW00} {Park, J. -H. and Won, S., Stability analysis for neutral delay-differential systems, \textit{
J. Franklin Inst.}, 2000, \textbf{337:} 1-9.}
\bibitem{QLWZ89} {Qin, Y. X., Liu, Y. Q., Wang, L., and Zheng, Z. X., {\it Stability of Dynamical Systems with Delays
(Second Edition)}, Science Press, Beijing, 1989.}

\bibitem{RL02} {Ren, H. S. and Li, H. Y., Explicit asymptotic stability criteria for neutral differential equations with
two delays, \textit{ Appl. Math. E-Notes}, 2002, \textbf{2:} 1-9.}

\bibitem{RZ99} {Ren, H. S. and Zheng, Z. X., Algebraic criterion for asymptotic stability of neutral equations of
the form $\dot{x}(t)+c\dot{x}(t-\tau)+ax(t)+bx(t-\tau)=0$,
\textit{ Acta. Math. Sinica}, 1999, \textbf{42(6):} 1077-1088.}

\bibitem{VN97} {Verriest, E. I. and Niculescu, S. I., Delay-independent stability of linear neutral systems:
A Riccati equation approach, in: {\it Stability and Control of Time-Delay Systems} (eds, L. Dugard and E. I. Verriest),
Springer-Verlag, London, 1997, 92-100.}

\bibitem{WLB02} {Wang, Z. D., Lam, J. and Burnham, K. J., Stability analysis and observer design for neutral delay systems,
\textit{ IEEE Trans. Automat. Contr.}, 2002, \textbf{AC-47:} 478-483.}


\bibitem{Wu94} {Wu, J. Y.,£¬ The position of boundary surface for a sort of trancendental equation with three
parameters, \textit{Acta. Math. Sinica}, 1994, \textbf{37:}
301-308.}

\bibitem{Yang99}  {Yang, L.,  Recent advances on determining the number of real roots of parametric polynomials,
\textit{J. Symbolic Computation,} 1999, \textbf{28:} 225-242.}

\bibitem{YHZ96}  {Yang, L., Hou, X. R. and Zeng, Z. B., A complete discrimination system for polynomials,
\textit{Science in China,} 1996, \textbf{E-39:} 628-646.}


\bibitem{YZH96} {Yang, L., Zhang, J. and Hou, X., {\it Nonliear Algebraic Equations and Machine Proving},
Shanghai Science and Education Press, Shanghai, 1996.}

\bibitem{YW00} {Yu, W. S. and Wang, L., Algebraic criteria for delay-independent stability of differential
  system $\dot{x}(t)=A_{0}x(t)+\sum^{N}_{k=1}A_{K}x(t-k\tau)$,
{\it Chinese Control Conference (CCC2000)}, Hong Kong, 2000,
435-439.}

\end{thebibliography}
\end{document}